\documentclass[12pt,reqno,openbib,runningheads,a4paper]{amsart}%
\usepackage{graphicx}
\usepackage{amssymb}
\usepackage{amsfonts}
\usepackage{amsmath}
\usepackage{graphicx}
\usepackage{lineno}
\usepackage[authoryear]{natbib}
\usepackage[dvipsnames]{xcolor}
\usepackage{epsf}
\usepackage{lscape}
\usepackage[legalpaper,bookmarks=true,colorlinks=true,linkcolor=blue,citecolor=blue]%
{hyperref}%
\providecommand{\U}[1]{\protect\rule{.1in}{.1in}}
\providecommand{\U}[1]{\protect\rule{.1in}{.1in}}
\newtheorem{theorem}{Theorem}[section]
\newtheorem{corollary}{Corollary}[section]

\newtheorem{lemma}{Lemma}[section]

\newtheorem{remark}{Remark}[section]

\makeatletter
\renewcommand{\@biblabel}[1]{}
\makeatother
\begin{document}

\begin{center}
{\Large \textbf{Gaussian approximation to the extreme value index estimator of
a heavy-tailed distribution}}

{\Large \textbf{under random censoring}}\medskip\medskip

{\large Brahim Brahimi, Djamel Meraghni, Abdelhakim Necir}$^{\ast}$\medskip

{\small \textit{Laboratory of Applied Mathematics, Mohamed Khider University,
Biskra, Algeria}}\medskip\medskip%
\[
\]

\end{center}

\noindent\textbf{Abstract}\medskip

\noindent We make use of the empirical process theory to approximate the
adapted Hill estimator, for censored data, in terms of Gaussian processes.
Then, we derive its asymptotic normality, only under the usual second-order
condition of regular variation. Our methodology allows to relax the
assumptions, made in \cite{EnFG08}, on the heavy-tailed distribution functions
and the sample fraction of upper order statistics.\medskip\medskip

\medskip


\vfill

\noindent{\small $^{\text{*}}$Corresponding author:
\texttt{necirabdelhakim@yahoo.fr} \newline\noindent\textit{E-mail
addresses:}\newline\texttt{brah.brahim@gmail.com} (B.~Brahimi)\newline%
\texttt{djmeraghni@yahoo.com} (D.~Meraghni)}

\section{\textbf{Introduction\label{sec1}}}

\noindent For $n\geq1,$ let $X_{1},X_{2},...,X_{n}$ be $n$ independent copies
of a non-negative continuous random variable (rv) $X,$ defined over some
probability space $\left(  \Omega,\mathcal{A},\mathbb{P}\right)  ,$ with
cumulative distribution function (cdf) $F.\ $We assume that the distribution
tail $1-F$ is regularly varying at infinity,\ with index $\left(
-1/\gamma_{1}\right)  ,$ notation: $1-F\in\mathcal{RV}_{\left(  -1/\gamma
_{1}\right)  }.$ That is%
\begin{equation}
\lim_{t\rightarrow\infty}\frac{1-F\left(  tx\right)  }{1-F\left(  t\right)
}=x^{-1/\gamma_{1}},\text{ for any }x>0, \label{first-condition}%
\end{equation}
where $\gamma_{1}>0,$ called shape parameter or tail index or extreme value
index (EVI), is a very crucial parameter in the analysis of extremes. It
governs the thickness of the distribution right tail: the heavier the tail,
the larger $\gamma_{1}.$ Its estimation has got a great deal of interest for
complete samples, as one might see in the textbook of \cite{BeGeS04}.\ In this
paper, we focus on the most celebrated (consistent and asymptotically normal)
estimator of $\gamma_{1},$ that was proposed by \cite{Hill75}:%
\[
\widehat{\gamma}_{1}^{H}=\widehat{\gamma}_{1}^{H}\left(  k\right)  :=\frac
{1}{k}%
%TCIMACRO{\dsum \limits_{i=1}^{k}}%
%BeginExpansion
{\displaystyle\sum\limits_{i=1}^{k}}
%EndExpansion
\log X_{n-i+1:n}-\log X_{n-k:n},
\]
where $X_{1:n}\leq...\leq X_{n:n}$ are the order statistics pertaining to the
sample $\left(  X_{1},...,X_{n}\right)  $ and $k=k_{n}$ is an integer sequence
satisfying%
\begin{equation}
1<k<n,\text{ }k\rightarrow\infty\text{ and }k/n\rightarrow0\text{ as
}n\rightarrow\infty. \label{K}%
\end{equation}
In the analysis of lifetime, reliability or insurance data, the observations
are usually randomly censored. In other words, in many real situations the
variable of interest $X$ is not always available. An appropriate way to model
this matter, is to introduce a non-negative continuous rv $Y,$ called
censoring rv, independent of $X$ and then to consider the rv $Z:=\min\left(
X,Y\right)  $ and the indicator variable $\delta:=\mathbf{1}\left(  X\leq
Y\right)  ,$ which determines whether or not $X$ has been observed. The cdf's
of $Y$ and $Z$ will be denoted by $G$ and $H$ respectively. The analysis of
extreme values of randomly censored data is a new research topic to which
\cite{ReTo97} made a very brief reference, in Section 6.1, as a first step but
with no asymptotic results. Considering Hall's model, \cite{BeGDFF07} proposed
estimators for the EVI and high quantiles and discussed their asymptotic
properties, when the data are censored by a deterministic threshold. More
recently, \cite{EnFG08} adapted various EVI estimators to the case where data
are censored, by a random threshold, and proposed a unified method to
establish their asymptotic normality. The obtained estimators are then used in
the estimation of extreme quantiles under random censorship. Moreover, they
applied their results on the Australian aids survival data available in the
MASS-package of the R software. \cite{GN11} also made a contribution to this
field by providing a detailed simulation study and applying the estimation
procedures on some survival data sets. In the same context, \cite{Worms}
presented a new approach, based on Kaplan-Meier integration, to define an
estimator for positive tail index and prove its consistency.\medskip

\noindent We start by a reminder of the definition of the adapted Hill
estimator, of the tail index $\gamma_{1},$ under random censorship.\ The tail
of the censoring distribution is assumed to be regularly varying too, that is
$1-G\in\mathcal{RV}_{\left(  -1/\gamma_{2}\right)  },$ for some $\gamma
_{2}>0.$\ By virtue of the independence of $X$ and $Y,$ we have $1-H\left(
x\right)  =\left(  1-F\left(  x\right)  \right)  \left(  1-G\left(  x\right)
\right)  $ and therefore $1-H\in\mathcal{RV}_{\left(  -1/\gamma\right)  }%
,$\ with $\gamma:=\gamma_{1}\gamma_{2}/\left(  \gamma_{1}+\gamma_{2}\right)
.$ Let $\left\{  \left(  Z_{i},\delta_{i}\right)  ,\text{ }1\leq i\leq
n\right\}  $ be a sample from the couple of rv's $\left(  Z,\delta\right)  $
and $Z_{1:n}\leq...\leq Z_{n:n}$\ the order statistics pertaining to $\left(
Z_{1},...,Z_{n}\right)  .$ In the sequel, the functions%
\[
H^{\left(  j\right)  }\left(  v\right)  :=\mathbb{P}\left(  Z\leq v,\text{
}\delta=j\right)  ,\text{ }j=0,1,
\]
will play a prominent role. If we denote the concomitant of the $i$th order
statistic by $\delta_{\left[  i:n\right]  }$ (i.e. $\delta_{\left[
i:n\right]  }=\delta_{j}$ if $Z_{i:n}=Z_{j}),$\ then the adapted Hill
estimator of $\gamma_{1}$ is defined by $\widehat{\gamma}_{1}:=\widehat
{\gamma}^{H}/\widehat{p},$ where $\widehat{\gamma}^{H}$ represents Hill's
estimator of $\gamma$ and $\widehat{p}:=k^{-1}\sum_{i=1}^{k}\delta_{\left[
n-i+1:n\right]  }$ estimates $p:=\gamma/\gamma_{1}.$ \cite{EnFG08} established
the asymptotic normality of $\widehat{\gamma}_{1}$ by assuming that cdf's $F$
and $G$ are absolutely continuous and that the quantile function $U_{H}\left(
t\right)  :=H^{\leftarrow}\left(  1-1/t\right)  ,$ $t\geq1,$ (the
notation\textbf{ }$K^{\leftarrow}$\textbf{ }stands for the quantile function
pertaining to a cdf $K)$\ satisfies the second-order condition of regular
variation \citep[see][]{deHS96}, that we write as follows:%
\begin{equation}
\lim_{t\rightarrow\infty}\frac{U_{H}\left(  tx\right)  /U_{H}\left(  t\right)
-x^{\gamma}}{A^{\ast}\left(  t\right)  }=x^{\gamma}\frac{x^{\tau}-1}{\tau
},\text{ for all }x>0, \label{UH}%
\end{equation}
where $\tau<0$ is the second-order parameter and $A^{\ast}\left(  t\right)  $
is a function tending to zero and not changing sign near infinity. For a
discussion on the relationships between $\left(  \ref{UH}\right)  $ and other
representations of the second-order condition of regular variation, like that
used in \cite{EnFG08}, one refers to \cite{FGH07}. In addition, \cite{EnFG08}
made three conditions $\left[  \mathcal{H}_{1}\right]  $-$\left[
\mathcal{H}_{3}\right]  $ (equivalently stated below) on the sequence $k$ in
terms of $U_{H}$ and an auxiliary function%
\begin{equation}
\mathbf{p}\left(  z\right)  :=\frac{\overline{G}\left(  z\right)  f\left(
z\right)  }{\overline{G}\left(  z\right)  f\left(  z\right)  +\overline
{F}\left(  z\right)  g\left(  z\right)  }, \label{p(z)}%
\end{equation}
where $f$ and $g$ represent the respective densities of cdf's $F$ and $G$ and,
for any $\mathcal{S},$ $\overline{\mathcal{S}}\left(  x\right)  :=\mathcal{S}%
\left(  \infty\right)  -\mathcal{S}\left(  x\right)  .$

\begin{itemize}
\item[.] $\left[  \mathcal{H}_{1}\right]  :\sqrt{k}A^{\ast}\left(  n/k\right)
\rightarrow d_{1}<\infty.$

\item[.] $\left[  \mathcal{H}_{2}\right]  :\dfrac{1}{\sqrt{k}}%
%TCIMACRO{\dsum \limits_{i=1}^{k}}%
%BeginExpansion
{\displaystyle\sum\limits_{i=1}^{k}}
%EndExpansion
\left[  \mathbf{p}\left(  H^{\leftarrow}\left(  1-\dfrac{i}{n}\right)
\right)  -p\right]  \rightarrow d_{2}<\infty.$

\item[.] $\left[  \mathcal{H}_{3}\right]  :\sqrt{k}\omega_{n,k}\left(
C\right)  :=\sqrt{k}\sup\limits_{\left(  s,t\right)  \in\mathcal{D}_{n}\left(
C\right)  }\left\vert \mathbf{p}\left(  H^{\leftarrow}\left(  t\right)
\right)  -\mathbf{p}\left(  H^{\leftarrow}\left(  s\right)  \right)
\right\vert \rightarrow0,$ where $\mathcal{D}_{n}\left(  C\right)  :=\left\{
1-k/n\leq t<1;\left\vert t-s\right\vert \leq C\sqrt{k}/n,s<1\right\}  ,$
$C>0.$
\end{itemize}

\noindent The authors claim that, in the case $\gamma_{1},\gamma_{2}>0,$ the
function $\mathbf{p}\left(  z\right)  $ tends to $p,$ as $z\rightarrow\infty,$
which , from a theoretical point of view, does not seem obvious by only
assuming the regular variation of $\overline{F}$ and $\overline{G},$ hence
further assumptions (like for instance the regular variation of $f$ and $g)$
are needed. From Theorem 1.7.2 in \cite{Bin87} page 39, known as the
\textit{Monotone Density Theorem}, a necessary condition for $f$ and $g$ to be
regularly varying is that $f$ and $g$ are ultimately monotone. In this paper,
we give an alternative definition to $\mathbf{p}\left(  z\right)  $ which,
among other things, avoids us additional restrictions on the underlying
distributions$.$ More precisely, we set $\mathbf{p}^{\ast}\left(  z\right)
:=\overline{H}^{\left(  1\right)  }\left(  z\right)  /\overline{H}\left(
z\right)  $ which, from assertion (i) of Lemma \ref{Lem1}, tends to $p.$ With
this choice and the use of the first-order regular variation conditions, we
show the consistency of $\widehat{p},$ which was not addressed by
\cite{EnFG08} who only proved its asymptotic normality under conditions
$\left[  \mathcal{H}_{2}\right]  $ and $\left[  \mathcal{H}_{3}\right]  .$
This result will then lead to the consistency of $\widehat{\gamma}_{1}.$ On
the other hand, we adopt an approach based on the empirical processes to
actually provide two main results. First, we solve the problems of restriction
by considering a more general family of distributions that only are regularly
varying at infinity. Second, in the restricted class of distributions we relax
the conditions $\left[  \mathcal{H}_{1}\right]  $-$\left[  \mathcal{H}%
_{3}\right]  $ on the sample fraction $k$ and reduce their number. To be\ more
precise, we show that $\left[  \mathcal{H}_{1}\right]  $-$\left[
\mathcal{H}_{3}\right]  $ imply that $\sqrt{k}\left(  \dfrac{n}{k}\overline
{H}^{\left(  1\right)  }\left(  U_{H}\left(  n/k\right)  \right)  -p\right)  $
tends to $d_{2}<\infty$ and that the converse is not true (see the end of
Appendix).\vspace{0.2cm}

\noindent The rest of the paper is organized a follows.\ In Section
\ref{sec2}, we state our main results which consist Gaussian approximations to
the adapted Hill estimator $\widehat{\gamma}_{1}$ in addition to its
consistency. This contribution, which, to the best of our knowledge, is the
first of its kind, will be of great usefulness in a lot of applications of
extreme value theory under random censoring. The proofs are postponed to
Section \ref{sec3} and some results, that are instrumental to our needs, are
gathered in the Appendix.

\section{\textbf{Main results\label{sec2}}}

\noindent We notice that the asymptotic normality of extreme value theory
based estimators is achieved in the second-order framework. Thus, it seems
quite natural to suppose that cdf's $F$ and $G$ satisfy the well-known
second-order condition of regular variation. That is, we assume that there
exist a constant $\tau_{j}<0$ and a function $A_{j}^{\ast},$ $j=1,2,$ tending
to zero and not changing sign near infinity, such that for any $x>0$%
\begin{equation}%
\begin{tabular}
[c]{l}%
$\underset{t\rightarrow\infty}{\lim}\dfrac{U_{F}\left(  tx\right)
/U_{F}\left(  t\right)  -x^{\gamma_{1}}}{A_{1}^{\ast}\left(  t\right)
}=x^{\gamma_{1}}\dfrac{x^{\tau_{1}}-1}{\tau_{1}},\medskip$\\
$\underset{t\rightarrow\infty}{\lim}\dfrac{U_{G}\left(  tx\right)
/U_{G}\left(  t\right)  -x^{\gamma_{2}}}{A_{2}^{\ast}\left(  t\right)
}=x^{\gamma_{2}}\dfrac{x^{\tau_{2}}-1}{\tau_{2}}.$%
\end{tabular}
\ \ \label{second-order}%
\end{equation}

\noindent In addition to approximating $\sqrt{k}\left(  \widehat{\gamma}%
_{1}-\gamma_{1}\right)  ,$ we also provide asymptotic normal representations
of two other useful statistics, namely $\sqrt{k}\left(  \widehat{p}-p\right)
$ and $\sqrt{k}\left(  Z_{n-k:n}/h-1\right)  ,$ where%
\begin{equation}
h=h_{n}:=U_{H}\left(  n/k\right)  . \label{h}%
\end{equation}
For convenience, we set, for $t>1,$ $A_{j}\left(  t\right)  :=A_{j}^{\ast
}\left(  1/\overline{F}\left(  t\right)  \right)  ,$ $j=1,2.$

\begin{theorem}
\label{Theorem1}Assume that $\overline{F}\in\mathcal{RV}_{\left(
-1/\gamma_{1}\right)  }$ and $\overline{G}\in\mathcal{RV}_{\left(
-1/\gamma_{2}\right)  }$ and let $k=k_{n}$ be an integer sequence satisfying
$(\ref{K}).$ Then $\widehat{\gamma}_{1}\rightarrow\gamma_{1}$ in probability.
Assume further that the second-order conditions $(\ref{second-order})$ hold
and $\sqrt{k}A_{j}\left(  h\right)  =O\left(  1\right)  ,$ for $j=1,2,$ as
$n\rightarrow\infty.$ Then there exists a sequence of Brownian bridges
$\left\{  B_{n}\left(  s\right)  ;\text{ }0\leq s\leq1\right\}  $ such that%
\begin{equation}
\sqrt{k}\left(  \frac{Z_{n-k:n}}{h}-1\right)  =\gamma\sqrt{\frac{n}{k}%
}\mathbb{B}_{n}^{\ast}\left(  \frac{k}{n}\right)  +o_{\mathbb{P}}\left(
1\right)  , \label{result1-1}%
\end{equation}%
\begin{align}
&  \sqrt{k}\left(  \widehat{p}-p\right) \label{result1-2}\\
&  =\sqrt{\frac{n}{k}}\left(  q\mathbb{B}_{n}\left(  \frac{k}{n}\right)
-p\widetilde{\mathbb{B}}_{n}\left(  \frac{k}{n}\right)  \right)  -pq\left(
\frac{\gamma_{1}^{-1}\sqrt{k}A_{1}\left(  h\right)  }{1-p\tau_{1}}%
-\frac{\gamma_{2}^{-1}\sqrt{k}A_{2}\left(  h\right)  }{1-q\tau_{2}}\right)
+o_{\mathbb{P}}\left(  1\right)  ,\nonumber
\end{align}
and%
\begin{equation}
\sqrt{k}\left(  \widehat{\gamma}_{1}-\gamma_{1}\right)  =\gamma_{1}\sqrt
{\frac{n}{k}}\int_{0}^{1}s^{-1}\mathbb{B}_{n}^{\ast}\left(  \frac{k}%
{n}s\right)  ds-\frac{\gamma_{1}}{p}\sqrt{\frac{n}{k}}\mathbb{B}_{n}\left(
\frac{k}{n}\right)  +\frac{\sqrt{k}A_{1}\left(  h\right)  }{1-p\tau_{1}%
}+o_{\mathbb{P}}\left(  1\right)  , \label{result1-3}%
\end{equation}
where, $\mathbb{B}_{n}\left(  s\right)  :=B_{n}\left(  \theta\right)
-B_{n}\left(  \theta-ps\right)  ,$ for $0\leq s<\theta/p,$ $\widetilde
{\mathbb{B}}_{n}\left(  s\right)  :=-B_{n}\left(  1-qs\right)  ,$ for $0\leq
s\leq1$ and $\mathbb{B}_{n}^{\ast}\left(  s\right)  :=\mathbb{B}_{n}\left(
s\right)  +\widetilde{\mathbb{B}}_{n}\left(  s\right)  ,$ for $0\leq
s<\theta/p,$ are sequences of centred Gaussian processes, with $\theta
:=H^{\left(  1\right)  }\left(  \infty\right)  $ and $q:=1-p.$
\end{theorem}

\begin{corollary}
\label{Cor1}Assume that the conditions of Theorem \ref{Theorem1} hold, assume
further that $\sqrt{k}A_{1}\left(  h\right)  \rightarrow\lambda_{1},$ as
$n\rightarrow\infty.$ Then%
\[
\sqrt{k}\left(  \widehat{\gamma}_{1}-\gamma_{1}\right)  \overset
{d}{\rightarrow}\mathcal{N}\left(  \dfrac{\lambda_{1}}{1-p\tau_{1}}%
,\dfrac{\gamma_{1}^{2}}{p}\right)  ,\text{ as }n\rightarrow\infty,
\]
where $\mathcal{N}\left(  m,d^{2}\right)  $ designates the normal distribution
with mean $m$ and variance $d^{2}.$
\end{corollary}

\begin{remark}
\label{Rq1}Given that $p=\gamma/\gamma_{1},$ the asymptotic variance above is
exactly the same as that obtained by \cite{EnFG08}.
\end{remark}

\begin{theorem}
\label{Theorem2}Let $F$ and $G$ be two absolutely continuous cdf's with
ultimately monotone densities. Assume that $\overline{F}\in\mathcal{RV}%
_{\left(  -1/\gamma_{1}\right)  }$ and $\overline{G}\in\mathcal{RV}_{\left(
-1/\gamma_{2}\right)  }$ and that $U_{H}$ satisfies the second-order condition
$\left(  \ref{UH}\right)  .$ Let $k=k_{n}$ be an integer sequence such that
$(\ref{K})$ holds and both $\ \sqrt{k}A^{\ast}\left(  n/k\right)  $ and
$\sqrt{k}\left(  \dfrac{n}{k}\overline{H}^{\left(  1\right)  }\left(
h\right)  -p\right)  $ are asymptotically bounded. Then%
\[
\sqrt{k}\left(  \frac{Z_{n-k:n}}{h}-1\right)  =\gamma\sqrt{\frac{n}{k}%
}\mathbb{B}_{n}^{\ast}\left(  \frac{k}{n}\right)  +o_{\mathbb{P}}\left(
1\right)  ,
\]%
\[
\sqrt{k}\left(  \widehat{p}-p\right)  =\sqrt{\frac{n}{k}}\left(
q\mathbb{B}_{n}\left(  \frac{k}{n}\right)  -p\widetilde{\mathbb{B}}_{n}\left(
\frac{k}{n}\right)  \right)  +\sqrt{k}\left(  \dfrac{n}{k}\overline
{H}^{\left(  1\right)  }\left(  h\right)  -p\right)  +o_{\mathbb{P}}\left(
1\right)  ,
\]
and%
\[
\sqrt{k}\left(  \widehat{\gamma}_{1}-\gamma_{1}\right)  =\gamma_{1}\sqrt
{\frac{n}{k}}\int_{0}^{1}s^{-1}\mathbb{B}_{n}^{\ast}\left(  \frac{k}%
{n}s\right)  ds-\frac{\gamma_{1}}{p}\sqrt{\frac{n}{k}}\mathbb{B}_{n}\left(
\frac{k}{n}\right)  +\sqrt{k}R_{n}+o_{\mathbb{P}}\left(  1\right)  ,
\]
where $R_{n}:=p^{-1}\left\{  \dfrac{\gamma}{1-\tau}A^{\ast}\left(  n/k\right)
-\gamma_{1}\left(  \dfrac{n}{k}\overline{H}^{\left(  1\right)  }\left(
h\right)  -p\right)  \right\}  .$
\end{theorem}

\begin{corollary}
\label{Cor2}Assume that the conditions of Theorem \ref{Theorem2} hold, assume
further that $\sqrt{k}A^{\ast}\left(  n/k\right)  \rightarrow d_{1}<\infty$
and $\sqrt{k}\left(  \dfrac{n}{k}\overline{H}^{\left(  1\right)  }\left(
h\right)  -p\right)  \rightarrow d_{2}<\infty,$ as $n\rightarrow\infty.$ Then%
\[
\sqrt{k}\left(  \widehat{\gamma}_{1}-\gamma_{1}\right)  \overset
{d}{\rightarrow}\mathcal{N}\left(  \dfrac{\gamma d_{1}}{p\left(
1-\tau\right)  }-\dfrac{\gamma_{1}d_{2}}{p},\dfrac{\gamma_{1}^{2}}{p}\right)
,\text{ as }n\rightarrow\infty.
\]

\end{corollary}

\begin{remark}
\label{Rq2}Note that the asymptotic bias above agrees with that obtained by
\cite{EnFG08}. For example, when $\gamma>-\rho,$ we combine relations $\left(
2.3\right)  ,$ $\left(  2.10\right)  $ and $\left(  3.4\right)  $ of
\cite{FGH07} to deduce that $A^{\ast}\left(  x\right)  =\frac{\gamma\left(
1-\tau\right)  }{\tau+\gamma\left(  1-\tau\right)  }b\left(  x\right)  ,$
where $b\left(  x\right)  $ is a function, defined in page $214$ of
\cite{EnFG08}, in terms of the convergence rate $a_{2}\left(  x\right)  $ of
the second-order condition, given in $\left(  9\right)  .$ Now, observe that
$\left[  \mathcal{H}_{1}\right]  $ is expressed as $\sqrt{k}b\left(
n/k\right)  \rightarrow\alpha_{1}$ and that the constant $d_{2}$ in $\left[
\mathcal{H}_{2}\right]  $ is denoted by $\alpha_{2},$ in the second
assumption, in \cite{EnFG08}. This leads, after substitution, to the same bias.
\end{remark}

\section{\textbf{Proofs\label{sec3}}}

\noindent We begin by a brief introduction on some uniform empirical processes
under random censoring. The empirical counterparts of $H^{\left(  j\right)  }$
$\left(  j=0,1\right)  $ and the pertaining empirical processes are
respectively defined, for $v\geq0,$ by%
\[
H_{n}^{\left(  j\right)  }\left(  v\right)  :=\frac{1}{n}\sum\limits_{i=1}%
^{n}\mathbf{1}\left(  Z_{i}\leq v,\delta_{i}=j\right)  \text{ and }\sqrt
{n}\left(  \overline{H}_{n}^{\left(  j\right)  }\left(  v\right)
-\overline{H}^{\left(  j\right)  }\left(  v\right)  \right)  ,\text{ }j=0,1.
\]
The latter may be represented, almost surely, by a uniform empirical process.
Indeed, let us define $U_{i}:=\delta_{i}H^{\left(  1\right)  }\left(
Z_{i}\right)  +\left(  1-\delta_{i}\right)  \left(  \theta+H^{\left(
0\right)  }\left(  Z_{i}\right)  \right)  ,$ $i=1,...,n.$ From \cite{EnKo92},
the rv's $U_{1},...,U_{n}$ are iid $(0,1)$-uniform. The empirical cdf and the
uniform empirical process based upon $U_{1},...,U_{n}$ are respectively
denoted by%
\[
\mathbb{U}_{n}\left(  s\right)  :=\frac{1}{n}\sum\limits_{i=1}^{n}%
\mathbf{1}\left(  U_{i}\leq s\right)  \text{ and }\alpha_{n}\left(  s\right)
:=\sqrt{n}\left(  \mathbb{U}_{n}\left(  s\right)  -s\right)  ,\text{ }0\leq
s\leq1.
\]
We have almost surely $H_{n}^{\left(  0\right)  }\left(  v\right)
=\mathbb{U}_{n}\left(  H^{\left(  0\right)  }\left(  v\right)  +\theta\right)
-\mathbb{U}_{n}\left(  \theta\right)  ,$ for $0<H^{\left(  0\right)  }\left(
v\right)  <1-\theta$ and $H_{n}^{\left(  1\right)  }\left(  v\right)
=\mathbb{U}_{n}\left(  H^{\left(  1\right)  }\left(  v\right)  \right)  ,$ for
$0<H^{\left(  1\right)  }\left(  v\right)  <\theta$ \citep[see][]{DeEn96}. It
is easy to verify that almost surely%
\begin{equation}
\sqrt{n}\left(  \overline{H}_{n}^{\left(  1\right)  }\left(  v\right)
-\overline{H}^{\left(  1\right)  }\left(  v\right)  \right)  =\alpha
_{n}\left(  \theta\right)  -\alpha_{n}\left(  \theta-\overline{H}^{\left(
1\right)  }\left(  v\right)  \right)  ,\text{ for }0<\overline{H}^{\left(
1\right)  }\left(  v\right)  <\theta, \label{rep-H1}%
\end{equation}
and%
\begin{equation}
\sqrt{n}\left(  \overline{H}_{n}^{\left(  0\right)  }\left(  v\right)
-\overline{H}^{\left(  0\right)  }\left(  v\right)  \right)  =-\alpha
_{n}\left(  1-\overline{H}^{\left(  0\right)  }\left(  v\right)  \right)
,\text{ for }0<\overline{H}^{\left(  0\right)  }\left(  v\right)  <1-\theta.
\label{rep-H0}%
\end{equation}
Our methodology strongly relies on the well-known Gaussian approximation,
given by \cite{CsCsHM86} by in Corollary 2.1, which says that on the
probability space $\left(  \Omega,\mathcal{A},\mathbb{P}\right)  ,$ there
exists a sequence of Brownian bridges $\left\{  B_{n}\left(  s\right)  ;\text{
}0\leq s\leq1\right\}  $ such that for every $0\leq\xi<1/4,$%
\begin{equation}
\sup_{1/n\leq s\leq1-1/n}\frac{n^{\xi}\left\vert \alpha_{n}\left(  s\right)
-B_{n}\left(  s\right)  \right\vert }{\left[  s\left(  1-s\right)  \right]
^{1/2-\xi}}=O_{\mathbb{P}}\left(  1\right)  . \label{approx}%
\end{equation}
For the increments $\alpha_{n}\left(  \theta\right)  -\alpha_{n}\left(
\theta-s\right)  ,$ we will need an approximation of the same type as $\left(
\ref{approx}\right)  .$ Following similar arguments, mutatis mutandis, as
those used to in the proof of assertions $\left(  2.2\right)  $ of Theorem 2.1
and $\left(  2.8\right)  $ of Theorem 2.2 in \cite{CsCsHM86}, we may show
that, for every $0<\theta<1$ and $0\leq\xi<1/4,$ we have%
\begin{equation}
\sup_{1/n\leq s\leq\theta}\frac{n^{\xi}\left\vert \left\{  \alpha_{n}\left(
\theta\right)  -\alpha_{n}\left(  \theta-s\right)  \right\}  -\left\{
B_{n}\left(  \theta\right)  -B_{n}\left(  \theta-s\right)  \right\}
\right\vert }{s^{1/2-\xi}}=O_{\mathbb{P}}\left(  1\right)  . \label{approx2}%
\end{equation}
The following processes will be crucial to our needs:%
\begin{equation}
\beta_{n}\left(  v\right)  :=\sqrt{\frac{n}{k}}\left\{  \alpha_{n}\left(
\theta\right)  -\alpha_{n}\left(  \theta-\overline{H}^{\left(  1\right)
}\left(  Z_{n-k:n}v\right)  \right)  \right\}  ,\text{ for }0<\overline
{H}^{\left(  1\right)  }\left(  v\right)  <\theta, \label{betan}%
\end{equation}
and%
\begin{equation}
\widetilde{\beta}_{n}\left(  v\right)  :=-\sqrt{\frac{n}{k}}\alpha_{n}\left(
1-\overline{H}^{\left(  0\right)  }\left(  Z_{n-k:n}v\right)  \right)  ,\text{
for }0<\overline{H}^{\left(  0\right)  }\left(  v\right)  <1-\theta.
\label{beta-tild}%
\end{equation}

\subsection{Proof of Theorem \ref{Theorem1}}

\noindent We start by showing the consistency of estimator $\widehat{p}.$
After observing that $\widehat{p}=\dfrac{n}{k}\overline{H}_{n}^{\left(
1\right)  }\left(  Z_{n-k:n}\right)  ,$ we consider the following
decomposition:%
\begin{align}
\widehat{p}-p  &  =\frac{n}{k}\left(  \overline{H}_{n}^{\left(  1\right)
}\left(  Z_{n-k:n}\right)  -\overline{H}^{\left(  1\right)  }\left(
Z_{n-k:n}\right)  \right) \label{phate-p}\\
&  +\frac{n}{k}\left(  \overline{H}^{\left(  1\right)  }\left(  Z_{n-k:n}%
\right)  -\overline{H}^{\left(  1\right)  }\left(  h\right)  \right)  +\left(
\frac{n}{k}\overline{H}^{\left(  1\right)  }\left(  h\right)  -p\right)
.\nonumber
\end{align}
Note that, from $\left(  \ref{betan}\right)  ,$ we have almost surely
\begin{equation}
\sqrt{k}\frac{n}{k}\left(  \overline{H}_{n}^{\left(  1\right)  }\left(
Z_{n-k:n}\right)  -\overline{H}^{\left(  1\right)  }\left(  Z_{n-k:n}\right)
\right)  =\beta_{n}\left(  1\right)  . \label{beta-p}%
\end{equation}
By using the Gaussian approximation $\left(  \ref{approx}\right)  ,$ we get
\[
\sqrt{k}\frac{n}{k}\left(  \overline{H}_{n}^{\left(  1\right)  }\left(
Z_{n-k:n}\right)  -\overline{H}^{\left(  1\right)  }\left(  Z_{n-k:n}\right)
\right)  =\sqrt{n/k}\mathbf{B}_{n}\left(  Z_{n-k:n}\right)  +o_{\mathbb{P}%
}\left(  1\right)  ,
\]
where%
\begin{equation}
\mathbf{B}_{n}\left(  v\right)  :=B_{n}\left(  \theta\right)  -B_{n}\left(
\theta-\overline{H}^{\left(  1\right)  }\left(  v\right)  \right)  ,\text{ for
}0<\overline{H}^{\left(  1\right)  }\left(  v\right)  <\theta.
\label{Bn-etoil}%
\end{equation}
Making use of Lemma \ref{Lem2}, we infer that $\beta_{n}\left(  1\right)
=\sqrt{n/k}\mathbb{B}_{n}\left(  k/n\right)  +o_{\mathbb{P}}\left(  1\right)
.$ It is easy to verify that$\sqrt{n/k}\mathbb{B}_{n}\left(  k/n\right)  $ is
asymptotically centred Gaussian rv with variance\ $p,$ it follows that the
first term in the right-hand side of $\left(  \ref{phate-p}\right)  $ is
$O_{\mathbb{P}}\left(  k^{-1/2}\right)  $ which tends to zero in probability.
The second term in the right-hand side of $\left(  \ref{phate-p}\right)  $ may
be written as%
\begin{equation}
\frac{n}{k}\left(  \overline{H}^{\left(  1\right)  }\left(  Z_{n-k:n}\right)
-\overline{H}^{\left(  1\right)  }\left(  h\right)  \right)  =\frac{n}%
{k}\overline{H}^{\left(  1\right)  }\left(  h\right)  \left(  \frac
{\overline{H}^{\left(  1\right)  }\left(  Z_{n-k:n}\right)  }{\overline
{H}^{\left(  1\right)  }\left(  h\right)  }-1\right)  . \label{phate-p2}%
\end{equation}
Since $\overline{H}$ is regularly varying at infinity with index $\left(
-1/\gamma\right)  ,$ then by letting $v=1$ in part (i) of Lemma \ref{Lem1}, it
is easy to verify that $\overline{H}^{\left(  1\right)  }$ is also regularly
varying with the same index. First, we note that, combining Corollary 2.2.2
with Potter's inequalities given in Proposition B.1.9 (5) in \cite{deHF06},
yields that $Z_{n-k:n}/h\rightarrow1$ in probability. Now, we use Potter's
inequalities \citep[see][Proposition B.1.9 (5)]{deHF06} with the fact that
$Z_{n-k:n}/h\rightarrow1$ in probability, to show that $\overline{H}^{\left(
1\right)  }\left(  Z_{n-k:n}\right)  /\overline{H}^{\left(  1\right)  }\left(
h\right)  -1\rightarrow0,$ in probability. If, in addition to $v=1,$ we take
$z=h$ in part (i) of Lemma \ref{Lem1}, then we get $n\overline{H}^{\left(
1\right)  }\left(  h\right)  /k\rightarrow p$ and therefore the second term in
the right-hand side of $\left(  \ref{phate-p}\right)  $ tends to zero in
probability as well. Finally, noting that the third term in the right-hand
side of $\left(  \ref{phate-p}\right)  $ clearly goes to zero, yields that
$\widehat{p}-p$ in probability as sought.\medskip

\noindent For result $\left(  \ref{result1-1}\right)  ,$ we write%
\[
\frac{\overline{H}\left(  tx\right)  }{\overline{H}\left(  t\right)
}-x^{-1/\gamma}=\frac{\overline{G}\left(  tx\right)  }{\overline{G}\left(
t\right)  }\frac{\frac{\overline{F}\left(  tx\right)  }{\overline{F}\left(
t\right)  }-x^{-1/\gamma_{1}}}{A_{1}\left(  t\right)  }A_{1}\left(  t\right)
+x^{-1/\gamma_{1}}\frac{\frac{\overline{G}\left(  tx\right)  }{\overline
{G}\left(  t\right)  }-x^{-1/\gamma_{2}}}{A_{2}\left(  t\right)  }A_{2}\left(
t\right)  .
\]
From Theorem 2.3.9 in \cite{deHF06}, page 48, the first and the second
conditions in $\left(  \ref{second-order}\right)  $ are respectively
equivalent to%
\[
\frac{\frac{\overline{F}\left(  tx\right)  }{\overline{F}\left(  t\right)
}-x^{-1/\gamma_{1}}}{A_{1}\left(  t\right)  }\rightarrow x^{-1/\gamma_{1}%
}\dfrac{x^{\tau_{1}/\gamma_{1}}-1}{\gamma_{1}\tau_{1}}\text{ and }\frac
{\frac{\overline{G}\left(  tx\right)  }{\overline{G}\left(  t\right)
}-x^{-1/\gamma_{2}}}{A_{2}\left(  t\right)  }\rightarrow x^{-1/\gamma_{2}%
}\dfrac{x^{\tau_{2}/\gamma_{2}}-1}{\gamma_{2}\tau_{2}},
\]
as $t\rightarrow\infty.$ This implies that%
\begin{equation}
\frac{\overline{H}\left(  tx\right)  }{\overline{H}\left(  t\right)
}-x^{-1/\gamma}\sim x^{-1/\gamma}\dfrac{x^{\tau_{1}/\gamma_{1}}-1}{\gamma
_{1}\tau_{1}}A_{1}\left(  t\right)  +x^{-1/\gamma}\dfrac{x^{\tau_{2}%
/\gamma_{2}}-1}{\gamma_{2}\tau_{2}}A_{2}\left(  t\right)  . \label{equiva}%
\end{equation}
In the sequel, for two sequences of rv's, we write $V_{n}^{\left(  1\right)
}=o_{\mathbb{P}}\left(  V_{n}^{\left(  2\right)  }\right)  ,$\ as
$n\rightarrow\infty,$ to say that $V_{n}^{\left(  1\right)  }/V_{n}^{\left(
2\right)  }\rightarrow0$\ in probability. Let now $x=x_{n}=Z_{n-k:n}/h$ and
$t=t_{n}=h.$ Since $x_{n}=1+o_{\mathbb{P}}\left(  1\right)  ,$ then
$\dfrac{x_{n}^{-\tau_{i}/\gamma_{i}}-1}{\tau_{i}\gamma_{i}},$ $i=1,2,$ tend to
zero in probability and therefore by using $\left(  \ref{equiva}\right)  ,$ we
get%
\begin{equation}
\left(  \frac{Z_{n-k:n}}{h}\right)  ^{-1/\gamma}=\frac{\overline{H}\left(
Z_{n-k:n}\right)  }{\overline{H}\left(  h\right)  }+o_{\mathbb{P}}\left(
A_{1}\left(  h\right)  +A_{2}\left(  h\right)  \right)  . \label{Zsurh}%
\end{equation}
But $\overline{H}\left(  h\right)  =k/n,$ then%
\[
\left(  \frac{Z_{n-k:n}}{h}\right)  ^{-1/\gamma}-1=\frac{\overline{H}\left(
Z_{n-k:n}\right)  }{k/n}-1+o_{\mathbb{P}}\left(  A_{1}\left(  h\right)
+A_{2}\left(  h\right)  \right)  .
\]
Applying the mean value theorem to the left-hand side yields%
\[
-\frac{1}{\gamma}\left(  \frac{Z_{n-k:n}}{h}-1\right)  c_{n}^{-1/\gamma
-1}=\frac{\overline{H}\left(  Z_{n-k:n}\right)  }{k/n}-1+o_{\mathbb{P}}\left(
A_{1}\left(  h\right)  +A_{2}\left(  h\right)  \right)  ,
\]
where $c_{n}$ is a sequence of rv's lying between $1$ and $Z_{n-k:n}/h,$
meaning that we have $c_{n}=1+o_{\mathbb{P}}\left(  1\right)  .$ It follows
that%
\[
\frac{Z_{n-k:n}}{h}-1=-\left(  1+o_{\mathbb{P}}\left(  1\right)  \right)
\gamma\left(  \frac{\overline{H}\left(  Z_{n-k:n}\right)  }{k/n}-1\right)
+o_{\mathbb{P}}\left(  A_{1}\left(  h\right)  +A_{2}\left(  h\right)  \right)
.
\]
By assumption, we have $\sqrt{k}A_{j}\left(  h\right)  =O\left(  1\right)  ,$
$j=1,2,$ then%
\[
\sqrt{k}\left(  \frac{Z_{n-k:n}}{h}-1\right)  =\left(  1+o_{\mathbb{P}}\left(
1\right)  \right)  \gamma\sqrt{k}\frac{n}{k}\left(  k/n-\overline{H}\left(
Z_{n-k:n}\right)  \right)  +o_{\mathbb{P}}\left(  1\right)  .
\]
Recall that $\overline{H}_{n}\left(  Z_{n-k:n}\right)  =k/n,$ then%
\[
\sqrt{k}\left(  \frac{Z_{n-k:n}}{h}-1\right)  =\left(  1+o_{\mathbb{P}}\left(
1\right)  \right)  \gamma\sqrt{k}\frac{n}{k}\left(  \overline{H}_{n}\left(
Z_{n-k:n}\right)  -\overline{H}\left(  Z_{n-k:n}\right)  \right)
+o_{\mathbb{P}}\left(  1\right)  .
\]
The right-hand side of the previous equation may be decomposed into%
\begin{align*}
&  \gamma\sqrt{k}\frac{n}{k}\left(  \overline{H}_{n}^{\left(  1\right)
}\left(  Z_{n-k:n}\right)  -\overline{H}^{\left(  1\right)  }\left(
Z_{n-k:n}\right)  \right) \\
&  \ \ \ \ \ \ \ \ \ \ +\gamma\sqrt{k}\frac{n}{k}\left(  \overline{H}%
_{n}^{\left(  0\right)  }\left(  Z_{n-k:n}\right)  -\overline{H}^{\left(
0\right)  }\left(  Z_{n-k:n}\right)  \right)  \left(  1+o_{\mathbb{P}}\left(
1\right)  \right)  +o_{\mathbb{P}}\left(  1\right)  .
\end{align*}
Using $\left(  \ref{betan}\right)  $ and $\left(  \ref{beta-tild}\right)  $
with $v=1,$ leads to%
\begin{equation}
\sqrt{k}\left(  \frac{Z_{n-k:n}}{h}-1\right)  =\gamma\left(  \beta_{n}\left(
1\right)  +\widetilde{\beta}_{n}\left(  1\right)  \right)  \left(
1+o_{\mathbb{P}}\left(  1\right)  \right)  +o_{\mathbb{P}}\left(  1\right)  ,
\label{Zh}%
\end{equation}
which, by the Gaussian representations $\left(  \ref{approx}\right)  $ and
$\left(  \ref{approx2}\right)  ,$ becomes%
\begin{equation}
\sqrt{k}\left(  \frac{Z_{n-k:n}}{h}-1\right)  =\gamma\sqrt{\frac{n}{k}%
}\mathbf{B}_{n}^{\ast}\left(  Z_{n-k:n}\right)  \left(  1+o_{\mathbb{P}%
}\left(  1\right)  \right)  +o_{\mathbb{P}}\left(  1\right)  , \label{Zn-k}%
\end{equation}
where%
\begin{equation}
\mathbf{B}_{n}^{\ast}\left(  v\right)  :=\mathbf{B}_{n}\left(  v\right)
-B_{n}\left(  1-\overline{H}^{\left(  0\right)  }\left(  v\right)  \right)
,\text{ for }0<\overline{H}^{\left(  0\right)  }\left(  v\right)  <1-\theta,
\label{Bn}%
\end{equation}
with $\mathbf{B}_{n}\left(  v\right)  $ defined in $\left(  \ref{Bn-etoil}%
\right)  .$ Finally, we use assertion $\left(  ii\right)  $ of Lemma
\ref{Lem2} to complete the proof of result $\left(  \ref{result1-1}\right)  .$
For result $\left(  \ref{result1-2}\right)  ,$ we multiply decomposition
$\left(  \ref{phate-p}\right)  $ by $\sqrt{k}$ and get%
\begin{align}
\sqrt{k}\left(  \widehat{p}-p\right)   &  =\sqrt{k}\frac{n}{k}\left(
\overline{H}_{n}^{\left(  1\right)  }\left(  Z_{n-k:n}\right)  -\overline
{H}^{\left(  1\right)  }\left(  Z_{n-k:n}\right)  \right)  \label{phate-p-bis}%
\\
&  +\sqrt{k}\frac{n}{k}\left(  \overline{H}^{\left(  1\right)  }\left(
Z_{n-k:n}\right)  -\overline{H}^{\left(  1\right)  }\left(  h\right)  \right)
+\sqrt{k}\left(  \frac{n}{k}\overline{H}^{\left(  1\right)  }\left(  h\right)
-p\right)  .\nonumber
\end{align}
Next, we represent the first two terms, of the right-side hand of the previous
equation, by $\beta_{n}$ and $\widetilde{\beta}_{n}$ and then we use
approximations $\left(  \ref{approx}\right)  $ and $\left(  \ref{approx2}%
\right)  .$ Recall that, from $\left(  \ref{beta-p}\right)  ,$ the first one
has already been shown to be equal to $\beta_{n}\left(  1\right)  $ almost
surely. For the second one, we use $\left(  \ref{phate-p2}\right)  $ and
assertion (ii) of Lemma \ref{Lem1} to have%
\[
\frac{n}{k}\left(  \overline{H}^{\left(  1\right)  }\left(  Z_{n-k:n}\right)
-\overline{H}^{\left(  1\right)  }\left(  h\right)  \right)  =p\left\{
\left(  \frac{Z_{n-k:n}}{h}\right)  ^{-1/\gamma}-1\right\}  +o_{\mathbb{P}%
}\left(  A_{1}\left(  h\right)  +A_{2}\left(  h\right)  \right)  .
\]
By applying the mean value theorem and using the fact that $Z_{n-k:n}%
/h=1+o_{\mathbb{P}}\left(  1\right)  ,$ we readily verify that $\left(
Z_{n-k:n}/h\right)  ^{-1/\gamma}-1=\left(  1+o_{\mathbb{P}}\left(  1\right)
\right)  \gamma^{-1}\left(  1-Z_{n-k:n}/h\right)  .$ Hence%
\begin{align}
&  \frac{n}{k}\left(  \overline{H}^{\left(  1\right)  }\left(  Z_{n-k:n}%
\right)  -\overline{H}^{\left(  1\right)  }\left(  h\right)  \right)
\label{sec-term}\\
&  =-\frac{p}{\gamma}\left(  \frac{Z_{n-k:n}}{h}-1\right)  \left(
1+o_{\mathbb{P}}\left(  1\right)  \right)  +o_{\mathbb{P}}\left(  A_{1}\left(
h\right)  +A_{2}\left(  h\right)  \right)  .\nonumber
\end{align}
Using the assumptions $\sqrt{k}A_{j}\left(  h\right)  =O\left(  1\right)  ,$
$j=1,2,$ and $\left(  \ref{Zh}\right)  ,$ we obtain%
\begin{equation}
\sqrt{k}\frac{n}{k}\left(  \overline{H}^{\left(  1\right)  }\left(
Z_{n-k:n}\right)  -\overline{H}^{\left(  1\right)  }\left(  h\right)  \right)
=-p\left(  \beta_{n}\left(  1\right)  +\widetilde{\beta}_{n}\left(  1\right)
\right)  +o_{\mathbb{P}}\left(  1\right)  . \label{Z}%
\end{equation}
Now, we use result (ii) of Lemma \ref{Lem1} with the fact that $\sqrt{k}%
A_{j}\left(  h\right)  =O\left(  1\right)  ,$ $j=1,2,$ to rewrite the third
term in the right-hand side of $\left(  \ref{phate-p-bis}\right)  $ into%
\begin{equation}
\sqrt{k}\left(  \frac{n}{k}\overline{H}^{\left(  1\right)  }\left(  h\right)
-p\right)  =\sqrt{k}b_{1}\left(  h\right)  +o\left(  1\right)  . \label{mup}%
\end{equation}
Substituting results $\left(  \ref{beta-p}\right)  ,$ $\left(  \ref{Z}\right)
$ and $\left(  \ref{mup}\right)  $ in decomposition $\left(  \ref{phate-p}%
\right)  ,$ yields%
\begin{equation}
\sqrt{k}\left(  \widehat{p}-p\right)  =\beta_{n}\left(  1\right)  -p\left(
\beta_{n}\left(  1\right)  +\widetilde{\beta}_{n}\left(  1\right)  \right)
+\sqrt{k}b_{1}\left(  h\right)  +o_{\mathbb{P}}\left(  1\right)  ,
\label{p-hate}%
\end{equation}
and the Gaussian approximations $\left(  \ref{approx}\right)  $ and $\left(
\ref{approx2}\right)  $ imply that%
\[
\sqrt{k}\left(  \widehat{p}-p\right)  =\sqrt{\frac{n}{k}}\left(
\mathbf{B}_{n}\left(  Z_{n-k:n}\right)  -p\mathbf{B}_{n}^{\ast}\left(
Z_{n-k:n}\right)  \right)  +\sqrt{k}b_{1}\left(  h\right)  +o_{\mathbb{P}%
}\left(  1\right)  .
\]
The final form of result $\left(  \ref{result1-2}\right)  $ is then obtained
by applying assertions $\left(  i\right)  $ and $\left(  ii\right)  $ of Lemma
\ref{Lem2}. Finally, we focus on $\left(  \ref{result1-3}\right)  ,$ which
represents the main result of Theorem \ref{Theorem1}. It is readily checked
that we have%
\begin{equation}
\sqrt{k}\left(  \widehat{\gamma}_{1}-\gamma_{1}\right)  =\frac{1}{\widehat{p}%
}\sqrt{k}\left(  \widehat{\gamma}^{H}-\gamma\right)  -\frac{\gamma_{1}%
}{\widehat{p}}\sqrt{k}\left(  \widehat{p}-p\right)  . \label{decomp}%
\end{equation}
Recall that one way to define Hill's estimator $\widehat{\gamma}^{H}$ is to
use the limit
\[
\gamma=\lim_{t\rightarrow\infty}\int_{t}^{\infty}v^{-1}\frac{\overline
{H}\left(  v\right)  }{\overline{H}(t)}dv.
\]
Then, by replacing $\overline{H}$ by $\overline{H}_{n}$ and letting
$t=Z_{n-k:n},$ we write%
\[
\widehat{\gamma}^{H}=\frac{n}{k}\int_{Z_{n-k:n}}^{\infty}v^{-1}\overline
{H}_{n}(v)dv.
\]
For details, see for instance, \cite{deHF06} page 69. Writing $\overline
{H}_{n}(v)$ as the sum of $\overline{H}_{n}^{\left(  0\right)  }(v)$ and
$\overline{H}_{n}^{\left(  1\right)  }(v),$ we decompose $\hat{\gamma}%
^{H}-\gamma$ into the sum of the following three terms%
\[
T_{n1}:=\frac{n}{k}\int_{Z_{n-k:n}}^{\infty}v^{-1}\left(  \overline{H}%
_{n}^{\left(  0\right)  }(v)-\overline{H}^{\left(  0\right)  }(v)+\overline
{H}_{n}^{\left(  1\right)  }(v)-\overline{H}^{\left(  1\right)  }(v)\right)
dv,
\]%
\[
T_{n2}:=\frac{n}{k}\int_{Z_{n-k:n}}^{h}v^{-1}\overline{H}\left(  v\right)
dv\text{ and }T_{n3}:=\frac{n}{k}\int_{h}^{\infty}v^{-1}\overline{H}\left(
v\right)  dv-\gamma.
\]
Making a change of variables in the first term $T_{n1}$ and using $\left(
\ref{rep-H1}\right)  ,$ $\left(  \ref{rep-H0}\right)  ,$ $\left(
\ref{betan}\right)  $ and $\left(  \ref{beta-tild}\right)  ,$ we get almost
surely%
\[
\sqrt{k}T_{n1}=\int_{1}^{\infty}v^{-1}\left(  \beta_{n}\left(  v\right)
+\widetilde{\beta}_{n}\left(  v\right)  \right)  dv.
\]
For the second term $T_{n2},$ we apply the mean value theorem to have%
\[
T_{n2}=\frac{n}{k}\frac{\overline{H}\left(  z_{n}^{\ast}\right)  }{z_{n}%
^{\ast}}\left(  h-Z_{n-k:n}\right)  ,
\]
where $z_{n}^{\ast}$ is a sequence of rv's lying between $Z_{n-k:n}$ and
$h.$\ It is obvious that we have $z_{n}^{\ast}=\left(  1+o_{\mathbb{P}}\left(
1\right)  \right)  h,$ this implies that $\overline{H}\left(  z_{n}^{\ast
}\right)  =\left(  1+o_{\mathbb{P}}\left(  1\right)  \right)  k/n.$ It follows
that the right-hand side of the previous equation is equal to $\left(
1+o_{\mathbb{P}}\left(  1\right)  \right)  \left(  1-Z_{n-k:n}/h\right)  .$
Hence, from $\left(  \ref{Zh}\right)  ,$ we have%
\[
\sqrt{k}T_{n2}=-\gamma\left(  \beta_{n}\left(  1\right)  +\widetilde{\beta
}_{n}\left(  1\right)  \right)  +o_{\mathbb{P}}\left(  1\right)  .
\]
Finally, for $T_{n3},$ we use the second-order conditions $\left(
\ref{second-order}\right)  $ to get%
\[
\sqrt{k}T_{n3}\sim p^{2}\frac{\sqrt{k}A_{1}\left(  h\right)  }{1-p\tau_{1}%
}+q^{2}\frac{\sqrt{k}A_{2}\left(  h\right)  }{1-q\tau_{2}}=:\sqrt{k}%
b_{2}\left(  h\right)  .
\]
Since by assumption $\sqrt{k}A_{j}\left(  h\right)  =O\left(  1\right)  ,$
$j=1,2,$ as $n\rightarrow\infty,$ then%
\begin{align}
&  \sqrt{k}\left(  \hat{\gamma}^{H}-\gamma\right) \label{rephill}\\
&  =\int_{1}^{\infty}v^{-1}\left(  \beta_{n}\left(  v\right)  +\widetilde
{\beta}_{n}\left(  v\right)  \right)  dv-\gamma\left(  \beta_{n}\left(
1\right)  +\widetilde{\beta}_{n}\left(  1\right)  \right)  +\sqrt{k}%
b_{2}\left(  h\right)  +o_{\mathbb{P}}\left(  1\right)  .\nonumber
\end{align}
Combining $\left(  \ref{p-hate}\right)  $ and $\left(  \ref{rephill}\right)  $
with $\left(  \ref{decomp}\right)  $ yields%
\begin{align*}
&  \sqrt{k}\left(  \widehat{\gamma}_{1}-\gamma_{1}\right) \\
&  =\frac{1}{p}\int_{1}^{\infty}v^{-1}\left(  \beta_{n}\left(  v\right)
+\widetilde{\beta}_{n}\left(  v\right)  \right)  dv-\frac{\gamma_{1}}{p}%
\beta_{n}\left(  1\right)  +\sqrt{k}b\left(  h\right)  +o_{\mathbb{P}}\left(
1\right)  ,
\end{align*}
where%
\begin{equation}
b\left(  h\right)  :=\frac{1}{p}b_{2}\left(  h\right)  -\frac{\gamma_{1}}%
{p}b_{1}\left(  h\right)  =\frac{A_{1}\left(  h\right)  }{1-p\tau_{1}}.
\label{b-h}%
\end{equation}
Once again, by using the Gaussian approximations $\left(  \ref{approx}\right)
$ and $\left(  \ref{approx2}\right)  ,$ we obtain%
\begin{align}
&  \sqrt{k}\left(  \widehat{\gamma}_{1}-\gamma_{1}\right)
\label{rep-hill-censor2}\\
&  =\frac{1}{p}\sqrt{\frac{n}{k}}\int_{1}^{\infty}v^{-1}\mathbf{B}_{n}^{\ast
}\left(  vZ_{n-k:n}\right)  dv-\frac{\gamma_{1}}{p}\sqrt{\frac{n}{k}%
}\mathbf{B}_{n}\left(  Z_{n-k:n}\right)  +\sqrt{k}b\left(  k\right)
+o_{\mathbb{P}}\left(  1\right)  ,\nonumber
\end{align}
which, by using assertions $\left(  i\right)  $ and $\left(  iii\right)  $ of
Lemma \ref{Lem2}, may be rewritten into%
\begin{align}
&  \sqrt{k}\left(  \widehat{\gamma}_{1}-\gamma_{1}\right)
\label{rep-hill-censor3}\\
&  =\frac{\gamma}{p}\sqrt{\frac{n}{k}}\int_{0}^{1}s^{-1}\mathbb{B}_{n}^{\ast
}\left(  \frac{k}{n}s\right)  ds-\frac{\gamma_{1}}{p}\sqrt{\frac{n}{k}%
}\mathbb{B}_{n}\left(  \frac{k}{n}\right)  +\sqrt{k}b\left(  h\right)
+o_{\mathbb{P}}\left(  1\right)  .\nonumber
\end{align}
Actually, applying the Gaussian approximations $\left(  \ref{approx}\right)  $
and $\left(  \ref{approx2}\right)  $ here needs some usual manipulation on the
upper bound of integration as recently done in \cite{BMNY-2013}, we omit
details. Replacing $\gamma$ by $p\gamma_{1}$ in the front of the integral and
using $\left(  \ref{b-h}\right)  $ achieve the proof of the third result of
Theorem \ref{Theorem1}.\hfill$\Box$

\subsection{Proof of Corollary \ref{Cor1}}

From result $\left(  \ref{result1-3}\right)  $ of Theorem \ref{Theorem1}, we
deduce that $\sqrt{k}\left(  \widehat{\gamma}_{1}-\gamma_{1}\right)  $ is
asymptotically Gaussian with mean $\dfrac{\lambda_{1}}{1-p\tau_{1}},$ and
variance%
\[
\gamma_{1}^{2}\lim_{n\rightarrow\infty}\mathbf{E}\left[  \sqrt{\frac{n}{k}%
}\int_{0}^{1}s^{-1}\mathbb{B}_{n}^{\ast}\left(  \frac{k}{n}s\right)
ds-\frac{1}{p}\sqrt{\frac{n}{k}}\mathbb{B}_{n}\left(  \frac{k}{n}\right)
\right]  ^{2}.
\]
We check that the processes $\mathbb{B}_{n}\left(  s\right)  ,$ $\widetilde
{\mathbb{B}}_{n}\left(  s\right)  $ and $\mathbb{B}_{n}^{\ast}\left(
s\right)  $ satisfy $p^{-1}\mathbf{E}\left[  \mathbb{B}_{n}\left(  s\right)
\mathbb{B}_{n}\left(  t\right)  \right]  =\min\left(  s,t\right)  -pst,$
$q^{-1}\mathbf{E}\left[  \widetilde{\mathbb{B}}_{n}\left(  s\right)
\widetilde{\mathbb{B}}_{n}\left(  t\right)  \right]  =\min\left(  s,t\right)
-qst$ and $p^{-1}\mathbf{E}\left[  \mathbb{B}_{n}\left(  s\right)
\mathbb{B}_{n}^{\ast}\left(  t\right)  \right]  =\mathbf{E}\left[
\mathbb{B}_{n}^{\ast}\left(  s\right)  \mathbb{B}_{n}^{\ast}\left(  t\right)
\right]  =\min\left(  s,t\right)  -st.$ Then, by elementary calculation (we
omit details), we obtain $\gamma_{1}^{2}/p$ for the asymptotic variance of
$\sqrt{k}\left(  \widehat{\gamma}_{1}-\gamma_{1}\right)  ,$ which, since
$p=\gamma/\gamma_{1},$ is equal to $\gamma_{1}^{3}/\gamma.$\hfill$\Box$

\subsection{Proof of Theorem \ref{Theorem2}}

Recall that%
\[
\sqrt{k}\left(  \widehat{\gamma}_{1}-\gamma_{1}\right)  =\frac{1}{\widehat{p}%
}\sqrt{k}\left(  \widehat{\gamma}^{H}-\gamma\right)  -\frac{\gamma_{1}%
}{\widehat{p}}\sqrt{k}\left(  \widehat{p}-p\right)  .
\]
For the first term, we use the Gaussian approximation $\left(  \ref{rephill}%
\right)  $ with a bias equal to $\sqrt{k}\dfrac{\gamma}{p\left(
1-\tau\right)  }A^{\ast}\left(  n/k\right)  $ under the assumption that
$\sqrt{k}A^{\ast}\left(  n/k\right)  $ is bounded. Now we consider $\sqrt
{k}\left(  \widehat{p}-p\right)  .$ The third term of the right-hand side of
$\left(  \ref{phate-p}\right)  $ is assumed to be asymptotically bounded while
the first one is already approximated by a Gaussian process by only being in
the first-order framework of regular variation. For the second term, denoted
by $\Omega_{n},$ we apply the mean value theorem to get%
\[
\Omega_{n}=\frac{n}{k}h\left(  Z_{n-k:n}/h-1\right)  \left(  \overline
{H}^{\left(  1\right)  }\right)  ^{\prime}\left(  \eta_{n}\right)  ,
\]
with $\eta_{n}$ between $h$ and $Z_{n-k:n}.$ Since $\overline{F}%
\in\mathcal{RV}_{\left(  -1/\gamma_{1}\right)  }$ with $f$ ultimately
monotone, then $f\in\mathcal{RV}_{\left(  -1/\gamma_{1}-1\right)  }$ by
Theorem 1.7.2 in \cite{Bin87} page 39, it follows that $z\left(  \overline
{H}^{\left(  1\right)  }\right)  ^{\prime}\left(  z\right)  \sim-\gamma
_{1}^{-1}\overline{H}\left(  z\right)  .$ From Potter's inequalities and the
fact that $Z_{n-k:n}/h=1+o_{\mathbb{P}}\left(  1\right)  ,$ we infer that%
\[
h\left(  \overline{H}^{\left(  1\right)  }\right)  ^{\prime}\left(  \eta
_{n}\right)  =-\left(  1+o_{\mathbb{P}}\left(  1\right)  \right)  \gamma
_{1}^{-1}\overline{H}\left(  h\right)  .
\]
Observe that $k/n=\overline{H}\left(  h\right)  ,$ then%
\[
\Omega_{n}=-\left(  1+o_{\mathbb{P}}\left(  1\right)  \right)  \gamma_{1}%
^{-1}\left(  Z_{n-k:n}/h-1\right)  ,
\]
which coincides with the first term in the right-hand side of $\left(
\ref{sec-term}\right)  .$ From this point on, we proceed as in the proof of
Theorem \ref{Theorem1} to complete the proof.\hfill$\Box\medskip$

\noindent\textbf{Concluding notes\medskip}

\noindent The primary object of the present work consists in providing a
Gaussian limiting distribution for the estimator of the shape parameter of a
heavy-tailed distribution, under random censorship. Our methodology is based
on the approximation of the uniform empirical process by a sequence of
Brownian bridges. This approach relaxes the assumptions imposed in Theorem 1
of \cite{EnFG08} and reduces their number. It is noteworthy that, for $p=1$
(non censoring case), our main result perfectly agrees with the Gaussian
representation of the classical Hill estimator. The generalization of this
approximating procedure to the moment estimator (valid for any real-valued
EVI) adapted to random censorship in \cite{EnFG08}, will be considered in a
future work.

\section{\textbf{Appendix\label{sec4}}}

\begin{lemma}
\label{Lem1}(i) Assume that $\overline{F}\in\mathcal{RV}_{\left(
-1/\gamma_{1}\right)  }$ and $\overline{G}\in\mathcal{RV}_{\left(
-1/\gamma_{2}\right)  }$ and let $k:=k_{n}$ be an integer sequence satisfying
$(\ref{K}).$ Then%
\[
\underset{z\rightarrow\infty}{\lim}\sup_{v\geq1}\left\vert \overline
{H}^{\left(  1\right)  }\left(  zv\right)  /\overline{H}\left(  z\right)
-pv^{-1/\gamma}\right\vert =0.
\]

\end{lemma}

(ii) If further the second-order conditions of regular variation $\left(
\ref{second-order}\right)  $ hold, then we have, uniformly on $v\geq1,$%
\[
v^{1/\gamma}\frac{n}{k}\overline{H}^{\left(  1\right)  }\left(  hv\right)
=p+b_{1}\left(  v;h\right)  +o\left(  A_{1}\left(  h\right)  +A_{2}\left(
h\right)  \right)  ,\text{ as }n\rightarrow\infty,
\]
where%
\[
b_{1}\left(  v;h\right)  :=\frac{p}{\gamma_{1}\tau_{1}}\left(  \frac
{1-\tau_{1}}{1-p\tau_{1}}v^{\tau_{1}/\gamma_{1}}-1\right)  A_{1}\left(
h\right)  +\frac{p}{\gamma_{2}\tau_{2}}\left(  \frac{1}{1-q\tau_{2}}%
v^{\tau_{2}/\gamma_{2}}-1\right)  A_{2}\left(  h\right)  .
\]

\begin{proof}
It is easy to verify that, for $v\geq1,$ we have $H^{\left(  1\right)
}\left(  v\right)  =\int_{0}^{v}\overline{G}\left(  y\right)  dF\left(
y\right)  $ and%
\[
\frac{\overline{H}^{\left(  1\right)  }\left(  zv\right)  }{\overline
{H}\left(  z\right)  }=-\int_{v}^{\infty}\frac{\overline{G}\left(  hx\right)
}{\overline{G}\left(  h\right)  }d\frac{\overline{F}\left(  hx\right)
}{\overline{F}\left(  h\right)  }.
\]
To prove (i), we directly use Proposition B.1.10 in \cite{deHF06} page 369 to
the regularly varying functions $\overline{F}$ and $\overline{G}:$ for any
$0<\epsilon<1,$ there exists $n_{0}=n_{0}\left(  \epsilon\right)  ,$ such that
for all $n>n_{0}$ and $x\geq1,$%
\[
\left\vert \frac{\overline{F}\left(  hx\right)  }{\overline{F}\left(
h\right)  }-x^{-1/\gamma_{1}}\right\vert <\epsilon x^{-1/\gamma_{1}+\epsilon
}\text{ and }\left\vert \frac{\overline{G}\left(  hx\right)  }{\overline
{G}\left(  h\right)  }-x^{-1/\gamma_{2}}\right\vert <\epsilon x^{-1/\gamma
_{2}+\epsilon}.
\]
For (ii), we use the uniform inequalities (for the second-order regularly
varying functions) to both tails $\overline{F}$ and $\overline{G}%
\ $\citep[see, e.g., ][bottom of page 161]{deHF06}: for any $0<\epsilon<1,$
there exists $n_{0}=n_{0}\left(  \epsilon\right)  ,$ such that for all
$n>n_{0}$ and $x\geq1$%
\begin{equation}
\left\vert \frac{\frac{\overline{F}\left(  hx\right)  }{\overline{F}\left(
h\right)  }-x^{-1/\gamma_{1}}}{A_{1}\left(  h\right)  }-x^{-1/\gamma_{1}%
}\dfrac{x^{\tau_{1}/\gamma_{1}}-1}{\gamma_{1}\tau_{1}}\right\vert \leq\epsilon
x^{-1/\gamma_{1}+\tau_{1}/\gamma_{1}+\epsilon}. \label{inequa1}%
\end{equation}
and%
\begin{equation}
\left\vert \frac{\frac{\overline{G}\left(  hx\right)  }{\overline{G}\left(
h\right)  }-x^{-1/\gamma_{2}}}{A_{2}\left(  h\right)  }-x^{-1/\gamma_{2}%
}\dfrac{x^{\tau_{2}/\gamma_{2}}-1}{\gamma_{2}\tau_{2}}\right\vert \leq\epsilon
x^{-1/\gamma_{2}+\tau_{2}/\gamma_{2}+\epsilon}. \label{inequa2}%
\end{equation}
Therefore we omit details.
\end{proof}

\begin{lemma}
\label{Lem2}Under the assumptions of assertion (i) in Lemma \ref{Lem1}, we
have%
\[%
\begin{array}
[c]{l}%
(i)\text{ }\sqrt{\dfrac{n}{k}}\mathbf{B}_{n}\left(  vZ_{n-k:n}\right)
=\sqrt{\dfrac{n}{k}}\mathbb{B}_{n}\left(  \dfrac{k}{n}v^{-1/\gamma}\right)
+o_{\mathbb{P}}\left(  1\right)  ,\text{ for every }v\geq1,\medskip\\
\left(  ii\right)  \text{ }\sqrt{\dfrac{n}{k}}\mathbf{B}_{n}^{\ast}\left(
vZ_{n-k:n}\right)  =\sqrt{\dfrac{n}{k}}\mathbb{B}_{n}^{\ast}\left(  \dfrac
{k}{n}v^{-1/\gamma}\right)  +o_{\mathbb{P}}\left(  1\right)  ,\text{ for every
}v\geq1,\medskip\\
(iii)\text{ }\sqrt{\dfrac{n}{k}}%
%TCIMACRO{\dint _{1}^{\infty}}%
%BeginExpansion
{\displaystyle\int_{1}^{\infty}}
%EndExpansion
v^{-1}\mathbf{B}_{n}^{\ast}\left(  vZ_{n-k:n}\right)  dv=\gamma\sqrt{\dfrac
{n}{k}}%
%TCIMACRO{\dint _{0}^{1}}%
%BeginExpansion
{\displaystyle\int_{0}^{1}}
%EndExpansion
s^{-1}\mathbb{B}_{n}^{\ast}\left(  \dfrac{k}{n}s\right)  ds+o_{\mathbb{P}%
}\left(  1\right)  .
\end{array}
\]

\end{lemma}

\begin{proof}
First, note that we will omit the details for assertion $\left(  ii\right)  ,$
as it is shown by similar arguments than those used to prove $\left(
i\right)  .$ For this latter, we have to show that for a fixed $v\geq1,$%
\[
\sqrt{\frac{n}{k}}\left\{  B_{n}\left(  \theta-\overline{H}^{\left(  1\right)
}\left(  Z_{n-k:n}v\right)  \right)  -B_{n}\left(  \theta-p\frac{k}%
{n}v^{-1/\gamma}\right)  \right\}  =o_{\mathbb{P}}\left(  1\right)  .
\]
Indeed, let $\left\{  W_{n}\left(  t\right)  ;\text{ }0\leq t\leq1\right\}  $
be a sequence of Wiener processes defined on $\left(  \Omega,\mathcal{A}%
,\mathbb{P}\right)  $ so that
\begin{equation}
\left\{  B_{n}\left(  t\right)  ;\text{ }0\leq t\leq1\right\}  \overset{d}%
{=}\left\{  W_{n}\left(  t\right)  -tW_{n}\left(  1\right)  ;\ 0\leq
t\leq1\right\}  . \label{W}%
\end{equation}
Then, without loss of generality, we may write%
\begin{align*}
&  \sqrt{\frac{n}{k}}\left\{  B_{n}\left(  \theta-\overline{H}^{\left(
1\right)  }\left(  Z_{n-k:n}v\right)  \right)  -B_{n}\left(  \theta-p\frac
{k}{n}v^{-1/\gamma}\right)  \right\} \\
&  =\sqrt{\frac{n}{k}}\left\{  W_{n}\left(  \theta-\overline{H}^{\left(
1\right)  }\left(  Z_{n-k:n}v\right)  \right)  -W_{n}\left(  \theta-p\frac
{k}{n}v^{-1/\gamma}\right)  \right\} \\
&  \ \ \ \ \ \ \ \ \ \ \ \ \ \ \ -\sqrt{\frac{n}{k}}\left(  p\frac{k}%
{n}v^{-1/\gamma}-\overline{H}^{\left(  1\right)  }\left(  Z_{n-k:n}v\right)
\right)  W_{n}\left(  1\right)  .
\end{align*}
Next, we show that both terms in the right-hand side tend to zero (in
probability), as $n\rightarrow\infty.$ Let $v\geq1$ be fixed and recall that
$\overline{H}^{\left(  1\right)  }\left(  Z_{n-k:n}v\right)  =\left(
1+o_{\mathbb{P}}\left(  1\right)  \right)  pv^{-1/\gamma}k/n,$ then the result
follows for the second term. For the first one, we have, for given
$0<\eta,\epsilon<1$ (small enough) and for all large $n,$%
\[
\mathbb{P}\left(  \left\vert \frac{\overline{H}^{\left(  1\right)  }\left(
Z_{n-k:n}v\right)  }{v^{-1/\gamma}k/n}-p\right\vert >\eta^{2}\frac
{\epsilon^{2}}{4v^{1/\gamma}}\right)  <\epsilon/2.
\]
Observe now that%
\begin{align*}
&  \mathbb{P}\left(  \sqrt{\frac{n}{k}}\left\vert W_{n}\left(  \theta
-\overline{H}^{\left(  1\right)  }\left(  Z_{n-k:n}v\right)  \right)
-W_{n}\left(  \theta-p\frac{k}{n}v^{-1/\gamma}\right)  \right\vert
>\eta\right) \\
&  =\mathbb{P}\left(  \sqrt{\frac{n}{k}}\left\vert W_{n}\left(  \left\vert
\overline{H}^{\left(  1\right)  }\left(  Z_{n-k:n}v\right)  -p\frac{k}%
{n}v^{-1/\gamma}\right\vert \right)  \right\vert >\eta\right) \\
&  \leq\mathbb{P}\left(  \left\vert \frac{\overline{H}^{\left(  1\right)
}\left(  Z_{n-k:n}v\right)  }{v^{-1/\gamma}k/n}-p\right\vert >\eta^{2}%
\frac{\epsilon^{2}}{4v^{1/\gamma}}\right)  +\mathbb{P}\left(  \sup_{0\leq
t\leq\frac{\epsilon^{2}}{4}\frac{k}{n}}\left\vert W_{n}\left(  t\right)
\right\vert >\eta\sqrt{k/n}\right)  .
\end{align*}
It is clear that the first term, in the right-hand side of the inequality
above, tends to zero as $n\rightarrow\infty.$ On the other hand, since
$\left\{  W_{n}\left(  t\right)  ;\text{ }0\leq t\leq1\right\}  $ is a
martingale, then by using the classical Doob inequality, we have, for any
$u>0$ and $T>0,$%
\[
\mathbb{P}\left(  \sup_{0\leq t\leq T}\left\vert W_{n}\left(  t\right)
\right\vert >u\right)  \leq\frac{\mathbf{E}\left\vert W_{n}\left(  T\right)
\right\vert }{u}\leq\frac{\sqrt{T}}{u}.
\]
Letting $T=\eta^{2}\frac{\epsilon^{2}}{4}\frac{k}{n}$ and $u=\eta\sqrt{k/n},$
yields that%
\[
\mathbb{P}\left(  \sup_{0\leq t\leq\eta^{2}\frac{\epsilon^{2}}{4}\frac{k}{n}%
}\left\vert W_{n}\left(  t\right)  \right\vert >\eta\sqrt{k/n}\right)
\leq\epsilon/2.
\]
This completes the proof of assertion $\left(  i\right)  .$ To prove $\left(
iii\right)  ,$ let us first show that%
\begin{equation}
\sqrt{\dfrac{n}{k}}%
%TCIMACRO{\dint _{h}^{Z_{n-k:n}}}%
%BeginExpansion
{\displaystyle\int_{h}^{Z_{n-k:n}}}
%EndExpansion
v^{-1}\mathbf{B}_{n}^{\ast}\left(  v\right)  dv=o_{\mathbb{P}}\left(
1\right)  . \label{reste}%
\end{equation}
For fixed $0<\eta,\epsilon<1,$ we have%
\begin{align*}
&  \mathbb{P}\left(  \left\vert \sqrt{\dfrac{n}{k}}%
%TCIMACRO{\dint _{h}^{Z_{n-k:n}}}%
%BeginExpansion
{\displaystyle\int_{h}^{Z_{n-k:n}}}
%EndExpansion
v^{-1}\mathbf{B}_{n}^{\ast}\left(  v\right)  dv\right\vert >\eta\right) \\
&  \leq\mathbb{P}\left(  \left\vert \frac{Z_{n-k:n}}{h}-1\right\vert
>\epsilon\right)  +\mathbb{P}\left(  \sqrt{\dfrac{n}{k}}%
%TCIMACRO{\dint _{h}^{\left(  1+\epsilon\right)  h}}%
%BeginExpansion
{\displaystyle\int_{h}^{\left(  1+\epsilon\right)  h}}
%EndExpansion
v^{-1}\left\vert \mathbf{B}_{n}^{\ast}\left(  v\right)  \right\vert
dv>\eta\right)  .
\end{align*}
Since $Z_{n-k:n}=\left(  1+o_{\mathbb{P}}\left(  1\right)  \right)  h,$ then
the first term of the previous expression tends to zero as $n\rightarrow
\infty.$ On the other hand, it is easy to verify that $\mathbf{E}\left[
\mathbf{B}_{n}^{\ast}\left(  u\right)  \mathbf{B}_{n}^{\ast}\left(  v\right)
\right]  =\min\left(  \overline{H}\left(  u\right)  ,\overline{H}\left(
v\right)  \right)  -\overline{H}\left(  u\right)  \overline{H}\left(
v\right)  ,$ therefore $\mathbf{E}\left\vert \mathbf{B}_{n}^{\ast}\left(
v\right)  \right\vert \leq\sqrt{\overline{H}\left(  v\right)  }.$ It follows
that%
\[
\mathbf{E}\left\vert \sqrt{\dfrac{n}{k}}%
%TCIMACRO{\dint _{h}^{\left(  1+\epsilon\right)  h}}%
%BeginExpansion
{\displaystyle\int_{h}^{\left(  1+\epsilon\right)  h}}
%EndExpansion
v^{-1}\mathbf{B}_{n}^{\ast}\left(  v\right)  dv\right\vert \leq\sqrt{\dfrac
{n}{k}}%
%TCIMACRO{\dint _{h}^{\left(  1+\epsilon\right)  h}}%
%BeginExpansion
{\displaystyle\int_{h}^{\left(  1+\epsilon\right)  h}}
%EndExpansion
v^{-1}\sqrt{\overline{H}\left(  v\right)  }dv,
\]
which, by a change of variables, is equal to $%
%TCIMACRO{\dint _{1}^{1+\epsilon}}%
%BeginExpansion
{\displaystyle\int_{1}^{1+\epsilon}}
%EndExpansion
v^{-1}\sqrt{\dfrac{n}{k}\overline{H}\left(  hv\right)  }dv.$ Since
$\overline{H}\left(  hv\right)  $ is asymptotically equivalent to
$v^{-1/\gamma}k/n$ (uniformly in $v\geq1),$ then the latter integral is
$\sim2\gamma\left(  1-\left(  1+\epsilon\right)  ^{-1/2\gamma}\right)  $ which
tends to zero as $\epsilon\downarrow0.$ Consequently, $\left(  \ref{reste}%
\right)  $ is true and we have%
\begin{equation}
\sqrt{\dfrac{n}{k}}%
%TCIMACRO{\dint _{1}^{\infty}}%
%BeginExpansion
{\displaystyle\int_{1}^{\infty}}
%EndExpansion
v^{-1}\mathbf{B}_{n}^{\ast}\left(  Z_{n-k:n}v\right)  dv=\sqrt{\dfrac{n}{k}}%
%TCIMACRO{\dint _{1}^{\infty}}%
%BeginExpansion
{\displaystyle\int_{1}^{\infty}}
%EndExpansion
v^{-1}\mathbf{B}_{n}^{\ast}\left(  hv\right)  dv+o_{\mathbb{P}}\left(
1\right)  . \label{Bn*}%
\end{equation}
Then to have $(iii),$ we need to show that%
\[
\sqrt{\dfrac{n}{k}}%
%TCIMACRO{\dint _{1}^{\infty}}%
%BeginExpansion
{\displaystyle\int_{1}^{\infty}}
%EndExpansion
v^{-1}\left(  \mathbf{B}_{n}^{\ast}\left(  hv\right)  -\mathbb{B}_{n}^{\ast
}\left(  \dfrac{k}{n}v^{-1/\gamma}\right)  \right)  dv=o_{\mathbb{P}}\left(
1\right)  .
\]
Observe that%
\begin{align*}
&  \sqrt{\dfrac{n}{k}}%
%TCIMACRO{\dint _{1}^{\infty}}%
%BeginExpansion
{\displaystyle\int_{1}^{\infty}}
%EndExpansion
v^{-1}\left(  \mathbf{B}_{n}^{\ast}\left(  hv\right)  -\mathbb{B}_{n}^{\ast
}\left(  \dfrac{k}{n}v^{-1/\gamma}\right)  \right)  dv\\
&  =\sqrt{\dfrac{n}{k}}%
%TCIMACRO{\dint _{1}^{\infty}}%
%BeginExpansion
{\displaystyle\int_{1}^{\infty}}
%EndExpansion
v^{-1}\left(  B_{n}\left(  \theta-p\dfrac{k}{n}v^{-1/\gamma}\right)
-B_{n}\left(  \theta-\overline{H}^{\left(  1\right)  }\left(  hv\right)
\right)  \right)  dv\\
&  +\sqrt{\dfrac{n}{k}}%
%TCIMACRO{\dint _{1}^{\infty}}%
%BeginExpansion
{\displaystyle\int_{1}^{\infty}}
%EndExpansion
v^{-1}\left(  B_{n}\left(  1-q\dfrac{k}{n}v^{-1/\gamma}\right)  -B_{n}\left(
1-\overline{H}^{\left(  0\right)  }\left(  hv\right)  \right)  \right)
dv=:I_{n1}+I_{n2}.
\end{align*}
We only show that $I_{n1}=o_{\mathbb{P}}\left(  1\right)  ,$ as similar
arguments lead to the same result for $I_{n2}$ (due to the symmetry structure
of $I_{n1}$ and $I_{n2}).$ Making use of the representation $\left(
\ref{W}\right)  ,$ we write%
\begin{align*}
&  I_{n1}\overset{d}{=}\sqrt{\dfrac{n}{k}}%
%TCIMACRO{\dint _{1}^{\infty}}%
%BeginExpansion
{\displaystyle\int_{1}^{\infty}}
%EndExpansion
v^{-1}\left(  W_{n}\left(  \theta-p\dfrac{k}{n}v^{-1/\gamma}\right)
-W_{n}\left(  \theta-\overline{H}^{\left(  1\right)  }\left(  hv\right)
\right)  \right)  dv\\
&  \ \ \ \ \ \ \ \ \ \ \ \ \ \ \ \ \ \ \ +W_{n}\left(  1\right)  \sqrt
{\dfrac{n}{k}}%
%TCIMACRO{\dint _{1}^{\infty}}%
%BeginExpansion
{\displaystyle\int_{1}^{\infty}}
%EndExpansion
v^{-1}\left(  p\dfrac{k}{n}v^{-1/\gamma}-\overline{H}^{\left(  1\right)
}\left(  hv\right)  \right)  dv.
\end{align*}
The increment $\left\{  W_{n}\left(  \theta-p\dfrac{k}{n}v^{-1/\gamma}\right)
-W_{n}\left(  \theta-\overline{H}^{\left(  1\right)  }\left(  hv\right)
\right)  \right\}  $ of the Wiener process $W_{n}$ is $\mathcal{N}\left(
0,\overline{H}^{\left(  1\right)  }\left(  hv\right)  -p\dfrac{k}%
{n}v^{-1/\gamma}\right)  ,$ then the absolute value of the first term of the
right-hand side has an expectation which is less than or equal to
\[
\sqrt{\dfrac{n}{k}}%
%TCIMACRO{\dint _{1}^{\infty}}%
%BeginExpansion
{\displaystyle\int_{1}^{\infty}}
%EndExpansion
v^{-1}\sqrt{\overline{H}^{\left(  1\right)  }\left(  hv\right)  -p\dfrac{k}%
{n}v^{-1/\gamma}}dv.
\]
Likewise, the absolute value of the second term of right-hand side has an
expectation which is less than or equal to%
\[
\sqrt{\dfrac{n}{k}}%
%TCIMACRO{\dint _{1}^{\infty}}%
%BeginExpansion
{\displaystyle\int_{1}^{\infty}}
%EndExpansion
v^{-1}\left\vert \overline{H}^{\left(  1\right)  }\left(  hv\right)
-p\dfrac{k}{n}v^{-1/\gamma}\right\vert dv.
\]
Routine manipulations, using assertion (i) in Lemma \ref{Lem1}, on the last
two quantities complete the proof.
\end{proof}

\textbf{Proof of }$\left[  \mathcal{H}_{2}\right]  \mathbf{-}\left[
\mathcal{H}_{3}\right]  \mathbf{\Longrightarrow}\sqrt{k}\left\{  \dfrac{n}%
{k}\overline{H}^{\left(  1\right)  }\left(  h\right)  -p\right\}
\mathbf{\rightarrow d}_{2}$

First, we show that%
\begin{equation}
\frac{1}{\sqrt{k}}%
%TCIMACRO{\dsum \limits_{i=1}^{k}}%
%BeginExpansion
{\displaystyle\sum\limits_{i=1}^{k}}
%EndExpansion
\left[  \mathbf{p}\left(  H^{\leftarrow}\left(  1-\dfrac{i}{n}\right)
\right)  -p\right]  =\sqrt{k}\left\{  \frac{n}{k}\overline{H}^{\left(
1\right)  }\left(  h\right)  -p\right\}  +O\left(  \sqrt{k}\widetilde{\omega
}_{n,k}\right)  , \label{implic}%
\end{equation}
where $\widetilde{\omega}_{n,k}$ is a sequence such that $\sqrt{k}%
\widetilde{\omega}_{n,k}\rightarrow0,$ as $n\rightarrow\infty.$ For
convenience, we set $\varphi\left(  s\right)  :=\mathbf{p}\left(
H^{\leftarrow}\left(  1-s\right)  \right)  \mathbf{,}$ $0<s\leq1.$ From the
classical Riemann approximation, we have%
\begin{align*}
&  \left\vert \frac{1}{k}%
%TCIMACRO{\dsum \limits_{i=1}^{k}}%
%BeginExpansion
{\displaystyle\sum\limits_{i=1}^{k}}
%EndExpansion
\varphi\left(  u+i\frac{v-u}{k}\right)  -\frac{1}{v-u}\int_{u}^{v}%
\varphi\left(  s\right)  ds\right\vert \\
&  \leq\sup_{\left\{  u<x,y\leq v,\left\vert y-x\right\vert \leq\left(
v-u\right)  /k\right\}  }\left\vert \varphi\left(  y\right)  -\varphi\left(
x\right)  \right\vert ,\text{ for }u<v<\infty.
\end{align*}
By letting $u=0$ and $v=k/n,$ we get
\[
\left\vert \frac{1}{k}%
%TCIMACRO{\dsum \limits_{i=1}^{k}}%
%BeginExpansion
{\displaystyle\sum\limits_{i=1}^{k}}
%EndExpansion
\varphi\left(  \frac{i}{n}\right)  -\frac{n}{k}\int_{0}^{k/n}\varphi\left(
s\right)  ds\right\vert \leq\widetilde{\omega}_{n,k},
\]
where $\widetilde{\omega}_{n,k}:=\sup_{\left\{  0<x,y\leq k/n,\left\vert
y-x\right\vert \leq1/n\right\}  }\left\vert \varphi\left(  y\right)
-\varphi\left(  x\right)  \right\vert .$ That is, we have%
\[
\frac{1}{k}%
%TCIMACRO{\dsum \limits_{i=1}^{k}}%
%BeginExpansion
{\displaystyle\sum\limits_{i=1}^{k}}
%EndExpansion
\mathbf{p}\left(  H^{\leftarrow}\left(  1-\frac{i}{n}\right)  \right)
=\frac{n}{k}\int_{0}^{k/n}\mathbf{p}\left(  H^{\leftarrow}\left(  1-x\right)
\right)  dx+O\left(  \widetilde{\omega}_{n,k}\right)  .
\]
Recall $\left(  \ref{p(z)}\right)  $ and observe that $\mathbf{p}\left(
z\right)  =\left(  \overline{H}^{\left(  1\right)  }\right)  ^{\prime}\left(
z\right)  /\overline{H}^{\prime}\left(  z\right)  .$ Then we have%
\[
\int_{0}^{k/n}\mathbf{p}\left(  H^{\leftarrow}\left(  1-x\right)  \right)
dx=\int_{0}^{k/n}\frac{\left(  H^{\left(  1\right)  }\right)  ^{\prime}\left(
H^{\leftarrow}\left(  1-x\right)  \right)  }{H^{\prime}\left(  H^{\leftarrow
}\left(  1-x\right)  \right)  }dx,
\]
which, by the change of variables $z=H^{\leftarrow}\left(  1-x\right)  ,$
equals $\int_{h}^{\infty}\left(  H^{\left(  1\right)  }\right)  ^{\prime
}\left(  z\right)  dz=\overline{H}^{\left(  1\right)  }\left(  h\right)  .$
Thus $\left(  \ref{implic}\right)  $ readily follows. Next, we use the limit
$\sqrt{k}\omega_{n,k}\left(  C\right)  \rightarrow0$ of assumption $\left[
\mathcal{H}_{3}\right]  .$ Observe that, for $n$ sufficiently large and all
$C>0,$ we have $1/n<C\sqrt{k}/n,$ then it is easy to verify that
$\widetilde{\omega}_{n,k}<\omega_{n,k}\left(  C\right)  .$ Hence $\sqrt
{k}\widetilde{\omega}_{n,k}\rightarrow0,$ which yields that the second term of
the right-hand side of the previous equation tends to zero as well. Finally,
we use assumption $\left[  \mathcal{H}_{2}\right]  $ to conclude.\hfill$\Box$

\end{document}